\documentclass[twoside]{article}
\usepackage{amssymb,a4}
\usepackage{amsthm}
\def\pt{\partial}
\def\n{\mathbb{N}}
\def\z{\mathbb{Z}}
\def\c{\mathbb{C}}
\def\dim{\hbox{dim}}

\def\a{\alpha}

\def\F{\mathbb C}

\def\ll{\lambda}

\newfont{\df}{eufm10}

\def\ll{\lambda}

\def\vp{\varphi}

\def\de{\delta}
\def\dim{\hbox{\rm dim}\,}

\hoffset \voffset \oddsidemargin=60pt \evensidemargin=45pt
\title{\bf { HARISH-CHANDRA MODULES OVER THE TWISTED HEISENBERG-VIRASORO
ALGEBRA} \footnotetext{E-mail: liudong@hutc.zj.cn,
cpjiang@sjtu.edu.cn}}
\author{Dong Liu$^{a}$\quad and \quad Cuipo Jiang$^{b}$
\\ $^a$Department of Mathematics, Huzhou Teachers College\\ Zhejiang Huzhou, 313000, P.R.China
\\ $^b$Department of Mathematics, Shanghai Jiaotong University\\ Shanghai, 200240, P.R.China}
\date{ }
\begin{document}
\maketitle

\def\abstractname{ABSTRACT}
\begin{abstract}

In this paper, we classify all indecomposable Harish-Chandra
modules of the intermediate series over the twisted
Heisenberg-Virasoro algebra. Meanwhile, some bosonic modules are
also studied.


\end{abstract}

\baselineskip 14pt

\setlength{\parindent}{1.5em}

\setcounter{section}{0}

\newtheorem{theo}{Theorem}[section]
\newtheorem{defi}[theo]{Definition}
\newtheorem{lemm}[theo]{Lemma}
\newtheorem{coro}[theo]{Corollary}
\newtheorem{prop}[theo]{Proposition}

\section{Introduction}

\vskip -3mm

The twisted Heisenberg-Virasoro algebra has been first studied by
Arbarello et al. in Ref.\cite{ADKP}, where a connection is
established between the second cohomology of certain moduli spaces
of curves and the second cohomology of the Lie algebra of
differential operators of order at most one:
$$L_{HV}=\{f(t)\frac{d}{dt}+g(t)|f,g\in{\mathbb
C}[t,t^{-1}]\}.$$ Moreover, the twisted Heisenberg-Virasoro
algebra has some relations with the toroidal Lie algebras (see
Ref.\cite{FM}) and the $N = 2$ Neveu-Schwarz superalgebra, which
is one of the most important algebraic objects realized in
superstring theory (see Ref.\cite{KL}).

By definition, as a vector space over $\c$, the twisted
Heisenberg-Virasoro algebra ${\cal L}$ has a basis $\{L(m), I(m),
C_L, C_I, C_{LI}, m\in\z\}$, subject to the following relations:
\begin{eqnarray*}
&&[L(m), L(n)]=(n-m)L(m+n)+\de_{m+n, 0}{1\over 12}(m^3-m)C_L;\\
&&[I(m), I(n)]=n\de_{m+n, 0}C_{I};\\
&&[L(m), I(n)]=nI(m+n)+\de_{m+n, 0}(m^2-m)C_{LI};\\
&&[{\cal L}, C_L]=[{\cal L}, C_I]=[{\cal L}, C_{LI}]=0.
\end{eqnarray*}

Clearly the Heisenberg algebra $H=\c\{I(m), C_I\mid m\in\z\}$ and
the Virasoro algebra $Vir=\c\{L(m), C_L\mid m\in\z\}$ are
subalgebras of ${\cal L}$. Moreover, ${\cal
L}=\oplus_{m\in\z}{\mathcal L}_m$, where ${\cal L}_m=\c\{t^m\pt,
t^m\}\oplus\de_{m, 0}\c\{C_I, C_L, C_{LI}\}$, is a $\z$-graded Lie
algebra.

 Arbarello et al.
Ref.\cite{ADKP} also proved that when the central element of the
Heisenberg subalgebra acts in a non-zero way, an irreducible
highest weight ${\cal L}$-module is isomorphic to the tensor
product of an irreducible $Vir$-module and an irreducible
$H$-module. The structure of the irreducible representations for
${\cal L}$ at level zero was studied in Ref.\cite{B}. The
irreducible weight modules were first studied in Ref.\cite{C}, and
then classified in Ref.\cite{LZ}. Some relations between
representation theory of the twisted Heisenberg-Virasoro algebra
and that of toroidal Lie algebras was studied in Ref.\cite{FM},
Ref.\cite {FO}, and Ref.\cite{JJ}.

By definition, a Harish-Chandra module over the twisted
Heisenberg-Virasoro algebra ${\cal L}$ is an indecomposable weight
module $V=\oplus V^\ll$, where $V^\ll=\{v\in V\mid L(0)v=(\ll+a)
v\}$ (for some $a\in\c$) and each $V^\ll$ is finite-dimensional.
Moreover if each $V_\ll$ is at most one-dimensional, the
nontrivial Harish-Chandra module $V$ is called the intermediate
series. By $[L(0), L(n)]=nL(n)$ and $[L(0), I(n)]=nI(n)$ for all
$n\in\z^*$, we see that if $V^\ll\ne 0$, then
$V=\oplus_{n\in\z}V^{\ll+n}$ and ${\cal L}_mV^{\ll+n}\subset
V^{\ll+n+m}$. Hence $V$ has a structure of $\z$-graded module.

Some indecomposable modules of the intermediate series over the
Lie algebra of differential operators of order at most one, were
constructed in Ref.\cite{C}. This construction is based on one
kind of modules of the intermediate series over the Virasoro
algebra and also needs a strong condition $L(\pm1)v_k\ne 0$ for
all $k\in\z$, which are not including the other two important
kinds of $Vir$-modules. In this paper we obtain all Harish-Chandra
modules of the intermediate series over ${\cal L}$ based on all
three kinds of the $Vir$-modules of the intermediate series. We
obtain seven kinds of new indecomposable modules over the twisted
Heisenberg-Virasoro algebra. Especially, we construct one kind of
new indecomposable modules from one kind of decomposable modules
of the Virasoro algebra.

The above researches may be very helpful in the following
researches. It can be used to classify Harish-Chandra modules over
the $N = 2$ Neveu-Schwarz superalgebra, whose even part is just
the twisted Heisenberg-Virasoro algebra (some central elements are
not consider here, see Ref.\cite{FJS}). Moreover, It also may be
helpful in the classification of the grade left-symmetric algebra
structure over this algebra (see Ref.\cite{KCB} in the Virasoro
case).

Spinor representations for the affine Lie algebras were first
developed by Frenkel Ref.\cite{F} and Kac-Peterson Ref.\cite{KP}
independently. The idea is to use a Clifford algebra with
infinitely many generators to construct certain quadratic
elements, which, together with the identity element, span an
orthogonal affine Lie algebra. Thereafter, Feingold-Frenkel
Ref.\cite{FF} constructed the so-called fermionic or bosonic
representations for all classical affine Lie algebras by using
Clifford or Weyl algebras with infinitely many generators. Bosonic
representations for the Lie algebra of differential operators were
studied in many papers (Ref.\cite{KR}, \cite{W}, \cite{HL}, etc.).
In this paper we also construct the bosonic modules based on such
construction over the twisted Heisenberg-Virasoro algebra.

This paper is organized as follows. In Section 2, we construct a
bosonic representation for the twisted Heisenberg-Virasoro algebra
by the construction over the Lie algebra of differential
operators. In Section 3, we classify all Harish-Chandra modules
over the twisted Heisenberg-Virasoro algebra. Throughout this
paper, $\z, \,\z^*$ and $\c$ denote the sets of integers, non-zero
integers and complex numbers, respectively.

\section{Bosonic modules over the twisted Heisenberg-Virasoro algebra}

\vskip -3mm
\vspace{2mm}

Define ${\mathcal S}$ to be the unital associative algebra with
infinitely many generators: $a(n), \;a^*(n) \; (n\in\z)$ with
relations
$$[\,a(n), a(m)\,]=[\,a^*(n), a^*(m)\,]=0,\eqno(2.1)$$
$$[\,a(n), a^*(m)\,]=-\delta_{n+m, 0}.\eqno(2.2)$$

We define the normal ordering as follows.
$$:a(n)a^*(m):=\cases{a(n)a^*(m),  \quad  n\le m,\cr   a^*(m)a(n),\quad n>m,}\eqno(2.3)$$
for $\;n, \,m\in\z.$  Set
$$\theta(n)=\cases{1, \quad n>0,\cr  0, \quad n\le0.} \eqno(2.4)$$
Then
$$a(n)a^*(m)=:a(n)a^*(m):-\,\delta_{n+m, 0}\,\theta(n-m),\eqno(2.5) $$
and
$$[\,a(m)a^*(n), a(p)\,]=\delta_{n+p, 0}\,a(m), $$
$$[\,a(m)a^*(n), a^*(p)\,]=-\delta_{m+p, 0}\,a^*(n),\eqno(2.6) $$
for $\;m, \,n, \,p\in \z$.

Let ${\mathcal S}^+$ be the subalgebra generated by $a(n),
\,a^*(0), \,a^*(m)$ for $n, \,m>0$. Let ${\mathcal S}^-$ be the
subalgebra generated by $a(0), \,a(n), \,a^*(m)$ for $n, \,m<0$.
Those generators in ${\mathcal S}^+$ are called annihilation
operators while those in ${\mathcal S}^-$ are called creation
operators. Let $V$ be a simple ${\mathcal S}$-module containing an
element $v_0$, called a ``vacuum vector", and satisfying
$${\mathcal S}^+v_0=0.\eqno(2.7)$$
So all annihilation operators kill $v_0$ and
$$V={\mathcal S}^-v_0.\eqno(2.8)$$

Now we may construct a class of bosons on $V$. For any $m\in\z,
\,n\in\n$, set
$$f(m, n)=\sum_{i\in \z}(-i)^n:a(m-i)a^*(i):.\eqno(2.9)$$
Although $f(m, n)$ are infinite sums, they are well-defined as
operators on $V$. Indeed, for any vector $v\in V={\mathcal
S}^-v_0$, only finitely many terms in (2.9) can make a non-zero
contribution to $f(m, n)v$.

\vspace{2mm}

\begin{prop} (Ref.\cite{W}, \cite{HL})  For
$m_1, \,m_2\in\z, \,n_1, \,n_2\in\n$,  we have
\begin{eqnarray*}
[\,f(m_1, n_1), f(m_2, n_2)\,]
&=&\sum_{i=0}^{n_1}{n_1\choose i}\,m_2^i\,f(m_1+m_2, n_1+n_2-i)\\
&& - \; \sum_{j=0}^{n_2}{n_2\choose j}\,m_1^j\,f(m_1+m_2, n_1+n_2-j)\\
&& + \; \vp\bigr(f(m_1, n_1), f(m_2, n_2)\bigl),
\end{eqnarray*}
{\it where $\vp$ is given by}
$$\vp(f(m_1, n_1), f(m_2,  n_2))=\cases{0, \ \hbox{\it if}\ m_1=0,\cr
(-1)^{n_1+1}\delta_{m_1+m_2,
0}\sum_{i=1}^{m_1}(m_1-i)^{n_1}i^{n_2}, \ \hbox{\it if}\ m_1>0,\cr
(-1)^{n_1}\delta_{m_1+m_2,
0}\sum_{i=m_1}^{-1}(m_1-i)^{n_1}i^{n_2}, \ \hbox{\it if}\ m_1<0.
}$$
\end{prop}

Let $T=f(0, 0)$, then
$$[\,T, a(n)\,]=a(n),\quad  [\,T, a^*(n)\,]=-a^*(n), \eqno(2.13)$$ for all $n\in \z$.
For any $v=a(n_1)\cdots a(n_s)a^*(m_1)\cdots a^*(m_l)v_0\in V$,
noting that $Tv_0=0$, one has
$$Tv=(s-l)v.\eqno(2.14)$$

According to Proposition 2.1, we obtain

\begin{theo}  $V$ is a module for the Lie algebra
${\cal L}$ with central charge $C_I=2, \, C_{LI}=1/2, \, C_L=1$
under the action given by
$$\pi(L(m))=f(m, 1),\quad \pi(I(m))=f(m, 0),$$
$$\pi(C_L)=2\,\hbox{id},\quad \pi(C_{LI})={1\over 2}\,\hbox{id},\quad \pi(C_I)=1\,\hbox{id},$$
 for all $\;m\in \z$. Moreover,
$$V=\bigoplus_{k\in\z} V_k,$$ is completely reducible,
where $V_k$ is an eigenspace with eigenvalue $k$ of operator $T$,
and each component $V_k$ is irreducible as a ${\cal L}$-module.
\end{theo}
\vspace{2mm}

{\bf Proof.} Note that
$$\vp(f(m, 1), f(n,  1))=\cases{0, \
\hbox{\it if}\quad m=0\cr \delta_{m+n,
0}\sum_{i=1}^{m}(m-i)i=\de_{m+n,0}{m^3-m\over 6}, \ \hbox{\it
if}\quad m>0\cr -\delta_{m+n,
0}\sum_{i=m}^{-1}(m-i)i=\de_{m+n,0}{m^3-m\over 6}, \ \hbox{\it
if}\quad m<0 }.$$ So $$\vp(f(m, 1), f(n,
1))=\de_{m+n,0}{m^3-m\over 6}.$$

Similarly,  $$\vp(f(m, 1), f(n, 0))=\de_{m+n,0}{m(m-1)\over 2}.$$
$$\vp(f(m, 0), f(n, 0))=n\de_{m+n,0}.$$

Hence, Proposition 2.1 shows that $V$ is a ${\cal L}$-module with
central charge $C_L=2,\, C_{LI}=1/2,\, C_I=1$.

The second statement is essentially same as the proof of Theorem
2.3 in Ref.\cite{HL}.  $\rule[-.23ex\,]{1.0ex}{2.0ex}$

\section{Harish-Chandra modules of the intermediate series over the twisted Heienberg-Virasoro algebra}

Kaplansky-Santharoubane Ref.\cite {KS} in 1983 gave a
classification of $Vir$-modules of the intermediate series. There
are three families of {\it indecomposable modules of the
intermediate series} (i.e nontrivial indecomposable weight modules
with each weight space is at most one-dimensional) over the
Virasoro algebra. They are $Vir$-modules "without central charge".

(1) ${\mathcal A}_{a,\; b}=\sum_{i\in\z}\c v_i$:
$L(m)v_i=(a+i+bm)v_{m+i}$;

(2) ${\mathcal A}(a)=\sum_{i\in\z}\c v_i$: $L(m)v_i=(i+m)v_{m+i}$
if $i\ne 0$, $L(m)v_0=m(m+a)v_{m}$;

(3) ${\mathcal B}(a)=\sum_{i\in\z}\c v_i$:  $L(m)v_i=iv_{m+i}$ if
$i\ne -m$, $L(m)v_{-m}=-m(m+a)v_0$,  for some $a, b\in\c$.

When $a\notin\z$ or $b\ne0, 1$, it is well-known that the module
${\cal A}_{a,\; b}$ is simple.  In the opposite case the module
contains two simple subquotients namely the trivial module and
$\c[t, t^{-1}]/\c$. Denote the nontrivial simple subquotients of
${\mathcal A}_{a,\; b}$, ${\mathcal A}(a)$, ${\mathcal B}(a)$ by
${\mathcal A}_{a, \; b}'$, ${\mathcal A}(a)'$, ${\mathcal B}(a)'$
respectively. They are all Harish-Chandra modules of the
intermediate series over the Virasoro algebra.

 Furthermore, we have

\begin{lemm}Ref.\cite{KS}
Among the above modules, we have following module isomorphisms:

(i) ${\mathcal A}_{a,\; b}\cong {\mathcal A}_{a',\; b}$ if
$a-a'\in\z$,

(ii) ${\mathcal A}_{a,\; 0}\cong {\mathcal A}_{a',\; 1}$ for
$a\not\in\z$ with $a-a'\in\z$,

(iii) ${\mathcal A}(a)'\cong {\mathcal B}(b)'\cong {\cal A}_{0,
0}'$ for all $a,\; b\in\c$.

\end{lemm}

\begin{theo} Ref.\cite{OZ}
Let $V$ be a $\z$-graded $Vir$-module with $\dim V >1$ and
$\dim V_j\le 1$ for all $j\in\z$. Suppose there exists $a\in\c$
such that $L_0$ acts on $V_j$ as the scalar $a + j$. Then $V$ is
isomorphic to one of the following for appropriate $a, b$: (i)
${\mathcal A}_{a,\; b}'$, (ii) ${\mathcal A}_{0,0}'\oplus\c v_0$
as direct sum of $Vir$-modules, (iii) ${\mathcal A}(a)$, (iv)
${\mathcal B}(a)$.

\end{theo}

For $a, b, c, d\in \c$, we can similarly define some
$H_{Vir}$-modules with  $C_{LI}=C_I=C_{L}=0$ as follows:

(1) ${\mathcal A}_{a, b, c}$:  $L(n)v_t=(a+t+bn)v_{n+t}$,
$I(n)v_t=cv_{n+t}$;

(2) ${\mathcal A}(a,\; d)$:  $L(n)v_t=(t+n)v_{n+t}$ if $t\ne 0$,
$L(n)v_0=n(n+a)v_{n}$, $I(n)v_t=0$ if $t\ne0$, $I(n)v_0=dn v_n$;

(3) ${\mathcal  B}(a, \;d)$:  $L(n)v_t=tv_{n+t}$ if $t\ne -n$,
$L(n)v_{-n}=-n(n+a)v_0$, $I(n)v_t=0$ if $n+t\ne 0$, $I(n)v_{-n}=dn
v_0$.

(4) ${\mathcal  U}_{d}$:  $L(n)v_t=tv_{n+t}$, $I(n)v_t=0$ if
$n+t\ne0$, $I(n)v_{-n}=n dv_{0}$ ;

(5) ${\mathcal  V}_{d}$:  $L(n)v_t=(t+n)v_{n+t}$, $I(n)v_t=0$ if
$t\ne 0$, $I(n)v_0=n dv_{n}$;

(6) $\tilde{\mathcal U}_{d}:=\sum_{m\in  M}\c v_m$:
$L(n)v_t=tv_{n+t}$ if $n+t\ne 0$, $L(n)v_{-n}=0$, $I(n)v_t=0$ if
$n+t\ne 0$, $I(n)v_{-n}=dn v_0, \, d\in \c^*$.

(7) $\tilde{\mathcal V}_{c}:=\sum_{m\in  M}\c v_m$:
$L(n)v_t=tv_{n+t}$ if $n+t\ne 0$, $L(n)v_{-n}=0$, $I(n)v_t=0$ if
$t\ne 0$ or $n=t=0$; $I(n)v_{0}=c v_n, \, n\ne 0, c\in \c^*$.

It is easy to prove  that ${\mathcal A}_{a, b, c}$ is simple if
and only if $a\notin M$ or $b\notin\{0, 1\}$ or $c\neq 0$. For
$a\in\z$,  the module ${\cal A}_{a,\; 0,\; 0}$ has an irreducible
quotient ${\cal A}'_{a,\; 0,\; 0}={\cal A}_{a,\; 0,\;
0}/{\F}v_{-a}$ and  ${\cal A}_{a,\; 1,\; 0}$ contains
 an irreducible submodule ${\cal A}'_{a,\; 1,\; 0}=\oplus_{n\neq-a}{\mathbb
F}v_{n}$.

Clearly, the above modules are indecomposable and are not
isomorphic each other. Now we shall classify all such modules,
i.e., all indecomposable weight modules with all weight
multiplicities $\le 1$ over the generalized Heisenberg-Virasoro
algebra. First we have the following result.

\begin{theo}
Among the ${\cal L}$-modules ${\mathcal A}_{a, \;b,\; c}$ for $a,
b, c\in\c$, and their nontrivial subquotients, we have following
module isomorphisms:

(i) ${\mathcal A}_{a,\; b,\; c}\cong {\mathcal A}_{a',\; b',\; c}$
if $a-a'\in\z$,

(ii) ${\mathcal A}_{a,\; 0,\; c}\cong {\mathcal A}_{a',\; 1,\; c}$
for $a\notin\z$ with $a-a'\in\z$,

\end{theo}
\noindent{\bf Proof.} It is clear.\hfill
$\rule[-.23ex]{1.0ex}{2.0ex}$

Now we shall classify all Harish-Chandra modules over the twisted
Heisenberg-Virasoro algebra. First we have the following two
lemmas.

\begin{lemm}
 Let $V=\sum_{t\in\z}\F v_t$ be a $M$-graded ${\cal L}$-module
such that $L(n)v_t=(a+t+bn)v_{t+n}$ for all $n, t\in\z$.

(i) If $a+bn\ne 0$ for all $n\in\z$, and $b\ne 0, 1$. Then
$$I(m)v_t=cv_{m+t}$$ for some $c\in\F$.

(ii) If $a\not\in\z$ and $b=0$, then

(1) $I(m)v_t=cv_{m+t}$, in this case $V\cong {\mathcal A}_{a,
0, c}$\\
or

 (2) $I(n)v_t={t+a\over n+t+a}cv_{n+t}$ for all $m, t\in
M$, in this case $V\cong {\mathcal A}_{a, 1, c}$ by setting
$v_t'={1\over a+t}v_t$.
\end{lemm}
\noindent{\bf Proof.} (i) Since $V$ is a module of the Virasoro
algebra Vir$=\oplus_{m\in\z} L(m)$, it is clear that $C_{L}=0$
(cf. [SZ] for example). Suppose that $I(n)v_t=f(n, t)v_{n+t}$ for
all $n, t\in\z$. From $[L(n), I(m)]=m I(n+m)+\de_{n+m,
0}(n^2-n)C_{LI}$, we see that
$$f(m, t)L(n)v_{t+m}-f(m, n+t)(a+t+bn)v_{n+t+m}$$$$=
m f(n+m, t)v_{n+m+t}+\de_{n+m, 0}(n^2-n)C_{LI}v_{t}. \eqno(3.1)$$
In this case,
$$f(m, t)(a+t+m+bn)-f(m, n+t)(a+t+bn)=m f(n+m, t)+\de_{n+m, 0}(n^2-n)C_{LI}. \eqno(3.2)$$
Let $m=0$ in (3.2), then we have
$$f(0,t)=f(0,n+t). \eqno(3.3)$$
Let $t=0$ in (3.2), then
$$f(m, n)(a+bn)=
(a+m+bn)f(m, 0)-(m f(n+m, 0)+\de_{n+m, 0}(n^2-n)C_{LI}), \quad
n\ne 0.\eqno(3.4)$$ Setting $n=-2m$ in (3.4), we have
$$
(a-2m b)f(m,-2m)=(a+m-2m b)f(m,0)-m f(-m,0).\eqno(3.5)$$
 Let $n=t=-m$ in (3.2),
then
$$f(m, -m)(a-bm)-f(m, -2m)(a-(1+b)m)=m f(0, -m)+(m^2+m)C_{LI}. \eqno(3.6)$$
Setting $n=-m$ in (3.4), then $$(a-m b)f(m,-m)=(a+m-bm)f(m,0)-m
f(0,0)-(m^{2}+m)C_{LI}.\eqno(3.7)$$ From (3.6) and (3.7) we obtain
that
$$
-(a-m-bm)f(m,-2m) =2m
f(0,0)-(a+m-bm)f(m,0)+2(m^{2}+m)C_{LI}.\eqno(3.8)$$ Therefore from
(3.5) and (3.8)  we have
$$
(a+(1-3b)m)f(m,0)+(a-(1+b)m)f(-m,0)$$$$
=2(a-2bm)f(0,0)+2(a-2bm)(1+m)C_{LI}, \ m\neq 0.\eqno(3.9)$$
Replacing $m$ by $-m$ in (3.9), (3.9) becomes
$$(a+(1+b)m)f(m, 0)+(a-(1-3b)m)f(-m, 0)$$
$$=2(a+2bm)f(0, 0)+2(a+2bm)(1-m)C_{LI}, \quad m\ne 0.\eqno(3.10)$$
From (3.9) and (3.10), we obtain
$$2(m^{2}b-m^{2}b^{2})(f(m,0)-f(0,0))
=(a^{2}m-am^{2}+2m^{2}b-2m^{2}b^{2}-4m^{3}b^{2}+m^{2}ab)C_{LI}.\eqno(3.11)$$
Setting $t=m=-n$ in (3.2), we have
$$(a+2m-bm)f(m,m)=m
f(0,0)+(a+m-bm)f(m,0)+(m^{2}+m)C_{LI}.\eqno(3.12)$$ Let $m=n$ in
(3.4), then
$$
(a+bm)f(m,m)=(a+m+bm)f(m,0)-m f(2m,0).\eqno(3.13)$$ From (3.12)
and (3.13) we have
$$
(a+bm)(m f(0,0)+(m^{2}+m)C_{LI})-2(m^{2}+am)f(m,0)+m (a+2m-bm)
f(2m,0)=0.\eqno(3.14)$$ Applying (3.11) to $f(m,0)$ and $f(2m,0)$,
we have
$$
(6(b^{3}-b^{2})m^{5}+3a(b-b^{2})m^{4}-2b^{2}m^{4}-2ab^{2}m^{3}-{3\over
2}a^{2}bm^{3}-{3\over 2}a^{2}m^{2})C_{LI}=0.\eqno(3.15)$$ Since
the equation (3.15) has infinitely many solutions and $b(1-b)\neq
0$, considering the coefficient of $m^{5}$, we deduce that
$$
C_{LI}=0.$$ Therefore  $$ f(m, 0)=f(0, 0),\eqno(3.16)$$ for all
$m\in\z$. By (3.4) we obtain that $f(m, n)=f(0, 0)=c\in\F$ for all
$m, n\in\z$. It follows that $C_{I}=0$.

(ii) If $a\not\in\z$ and $b=0$, then (3.1)---(3.15) still  hold.
From (3.11) we know that
$$
(a^{2}m-am^{2})C_{LI}=0, \ m\neq 0,$$ i.e.,
$$am(a-m)C_{LI}=0, m\neq 0.$$
Since $a\notin M$, we have $C_{LI}=0$. Thus (3.2) and (3.4) become
$$
(a+t+m)f(m,t)-(a+t)f(m,n+t)=m f(m+n,t)\eqno(3.17)$$ and
$$
(a+m)f(m,0)-af(m,n)=m f(m+n,0).\eqno(3.18)$$ Setting $n=-m$ in
(3.18), we have
$$
(a+m)f(m,0)-af(m,-m)=m f(0,0).\eqno(3.19)$$
 Let $n=m=-t$ in (3.17), then
 $$
 af(m,-m)-(a-m)f(m,0)=m f(2m,-m).\eqno(3.20)$$
By (3.18), we have
$$
af(2m,-m)=(a+2m)f(2m,0)-2m f(m,0).\eqno(3.21)$$ From (3.19) and
(3.20), we have $$2m f(m,0)=m f(2m,-m)+m f(0,0).\eqno(3.22)$$
Therefore
$$
(a+2m)f(2m,0)-2(m+a)f(m,0)+af(0,0)=0.\eqno(3.23)$$ Furthermore, we
have
$$(km+a)f(km, 0)=k(m+a)f(m, 0)-(k-1)af(0,
0).\eqno(3.24)$$

By $[I(m), I(n)]=\de_{m+n, 0}n C_I$, we get
$$f(n, t)f(m, n+t)=f(m, t)f(n, m+t)+\de_{m+n, 0}n C_I.\eqno(3.25)$$
So
$$f(2m, 0)f(m, 2m)=f(m, 0)f(2m, m).\eqno(3.26)$$
By (3.18) we obtain
$$f(m, 2m)={a+m\over a}f(m, 0)-{m\over a}f(3m, 0).\eqno(3.27)$$
$$f(2m, m)={a+2m\over a}f(2m, 0)-{2m\over a}f(3m, 0).\eqno(3.28)$$
By (3.24), we have
$$(3m+a)f(3m, 0)=3(m+a)f(m, 0)-2af(0, 0).\eqno(3.29)$$
Combining (3.23) and (3.26)---(3.29) we obtain
$$f(m, 0)=f(0, 0), \quad \hbox{or}\quad f(m, 0)={a\over a+m}f(0, 0).$$
So by (3.18) we have
$$f(m, n)=f(0, 0), \quad \hbox{or}\quad f(m, n)={a+n\over a+m+n}f(0, 0)\eqno(3.30)$$
and the lemma is proved. \hfill $\rule[-.23ex]{1.0ex}{2.0ex}$

\noindent{\bf Remark.} 1. In the following  cases, we  can also
deduce that $C_{L}=C_{LI}=C_I=0$  as in Lemma 3.4. So in the
following discussions, we always assume that $C_{L}=C_{LI}=C_I=0$.

 2. In Ref.\cite{C}, there is also a similar result (see Proposition 5.1 in
Ref.\cite{C}) as Lemma 3.4 on the Lie algebra $L_{HV}$, the Lie
algebra of differential operators of order at most one. But it
needs a strong condition $L(\pm1)v_i\ne 0$ for all $i\in\z$ and
its proof is more complicated. The two kinds of modules for $b=0$
and $b=1$ in Proposition 5.1 of Ref.\cite{C} are isomorphic if
$\a\not\in\z$.

Now we prove the main theorem of this section.
\begin{theo}
Let $V=\sum_{m\in\z}\mathbb Fv_m$ be an indecomposable ${\mathcal
L}$-module such that ${\cal L}_m v_n\in \F v_{m+n}$ for all $m,
n\in\z$. Then $V$ is isomorphic to one of the modules ${\mathcal
A}_{a, b, c}$, ${\mathcal A}(a, c)$, ${\mathcal B}(a, c)$,
${\mathcal U}_d$, ${\mathcal V}_d$, $\tilde{\mathcal U}_{d}$,
$\tilde{\mathcal V}_{c}$ for appropriate $a, b, c, d\in{\mathbb
F}$.
\end{theo}
\noindent{\bf Proof.} First we suppose that $V$ is an
indecomposable module of the  Virasoro algebra Vir$=\oplus_{m\in
\z}L(m)$, then $V$ is isomorphic to one of the modules in Theorem
3.3 as a Vir-module. Suppose that $I(n)v_t=f(n, t)v_{n+t}$ for all
$t, n\in\z$.

\noindent{\bf Case I.} $L(n)v_t=(a+t+bn)v_{n+t}$ for all $n,
t\in\z$.

({\bf I}.1) Suppose that $a\not\in\z$ and $a+bn\ne 0$ for all
$n\in\z$. If $b\ne 0, 1$, then $f(m, n)=f(0, 0)$ for all $m,
n\in\z$ by Lemma 3.4.

If $b=0$, then $V$ is isomorphic to ${\mathcal A}_{a, 0, c}$ or
${\cal A}_{a, 1, c}$ by Lemma 3.4.

If $b=1$, then $V$ is also isomorphic to ${\mathcal A}_{a, 0, c}$
or ${\cal A}_{a, 1, c}$  since ${\cal A}_{a, 0}\cong {\cal A}_{a,
1}$ for $a\not\in\z$.

({\bf I}.2) $a\not\in\z$ and $a=bp$ for some $p\in
M\backslash\{0\}$. So $b\neq 0,1$. Therefore $f(m,0)=f(0,0)$ and
by (3.4)
$$f(m,n)=f(0,0) \ {\rm if} \  n+p\neq 0.\eqno(3.31)$$
 It follows from (3.2)
that $$f(0, -p)=f(0, 0).\eqno(3.32)$$ Setting $t=-p$, $n=-m$ and
$m\neq p$ in (3.2) and using (3.31)-(3.32), we have
$$f(m,-p)(a-p+m-bm)=f(0,0)(a-p+m-bm).$$
Since $a=bp$ and $m\neq p$, we have
 $f(m, -p)=f(0, 0)$
if $m\ne p$. Finally, setting $m=n=p, t=-p$ in (3.2) we have $f(p,
-p)=f(0, 0)$. Therefore $f(m, n)=f(0, 0)$ for all $m,n\in\z$.

({\bf I}.3) $a\in\z$. Since ${\cal A}_{a, b}\cong {\cal A}_{0,
b}$, so we can suppose that $L(n)v_t=(t+bn)v_{n+t}$ for all $n,
t\in\z$.

 ({\bf I}.3.1)\ $b\ne 0, 1$, then $a+bn\ne 0$ for all $n\ne 0$. So $f(m,
n)=f(0, 0)$ by Lemma 3.4. Therefore $V$ is isomorphic to ${\cal
A}_{0,b,c}$.

({\bf I}.3.2) \ $b=1$. In this case (3.3) still holds and (3.2)
and (3.4) becomes
$$f(m, t)(t+m+n)-f(m, n+t)(t+n)=m f(n+m, t)\eqno(3.33)$$
and
$$
n f(m,n)=-m f(n+m,0)+(n+m)f(m,0)\eqno(3.34)$$ Replacing $t$ by $n$
and letting $n=-m$ in (3.33) and using (3.3), we have $$n
f(m,n)=(n-m)f(m,n-m)+m f(0,0).\eqno(3.35)$$ Replacing $n$ by $n-m$
in (3.34), we have
$$(n-m)f(m,n-m)=-m f(n,0)+n f(m,0).\eqno(3.36)$$
From (3.34)-(3.36), we have
$$
f(n+m,0)=f(n,0)+f(m,0)-f(0,0).\eqno(3.37)$$ We deduce that there
exists $d\in\F$ such that
$$
f(m,0)=dm+f(0,0),\eqno(3.38)$$ for all $m\in\z$. Applying (3.37)
to (3.34) and using (3.38), we have
$$
f(m,n)=f(0,0), \ n\neq 0.$$ Therefore
  $$f(m, n)=\cases{f(0, 0),\quad
n\ne 0,\cr  dm+f(0,0), \quad\hfill n=0,\ m\ne 0.}$$ Since $[I(m),
I(n)]=m\delta_{m+n,0}C_{I}$, we have
$$
f(m,n+t)f(n,t)=f(n,m+t)f(m,t).$$ Let $ t=0$, $n\neq 0$, $m\neq 0$
and $m\neq n$, then
$$
f(0,0)(f(0,0)+dn)=f(0,0)(f(0,0)+dm).$$ Therefore if $f(0,0)\neq
0$, then $d=0$ and $V\cong{\mathcal  A}_{0,1,c}$, where
$c=f(0,0)\neq 0$. If $f(0,0)=0$, then $$f(m,n)=f(0,0)=0, \ n\neq
0.$$ Therefore
  $$f(m, n)=\cases{0,\quad
n\ne 0,\cr  dm, \quad\hfill n=0.}$$ and $V\cong{\mathcal V}_{d}$.

 ({\bf I}.3.3) \ $b=0$. (3.2) becomes
$$
(t+m)f(m,t)-tf(m,n+t)=m f(m+n,t).\eqno(3.39)$$ Let $t=0$ and
$n=-m$ in (3.39), then
$$f(m,0)=f(0,0).\eqno(3.40)$$
Setting $t=m, n=-m$ in (3.39) and using (3.40), we have
$$
f(m,m)=0.\eqno(3.41)$$ Let $n+t=-m$ in (3.39), then$$ tf(m,-m)+m
f(-t,t)=(t+m)f(m,t).\eqno(3.42)$$ Therefore
$$f(m,-m)+f(-m,m)=2f(0,0).\eqno(3.43)$$
Setting $t=-m-n$ in (3.39), we have
$$
(-n)f(m,-m-n)+(m+n)f(m,-m)=m f(m+n,-m-n).\eqno(3.44)$$ Replace $m$
by $n'$ and $n$ by $m'$ respectively, then$$
(-m')f(n',-m'-n')+(m'+n')f(n',-n')=n'f(m'+n',-m'-n').\eqno(3.45)$$
Let $t=-m$ and replace $n$ by $-n$, then
$$
f(m,-n-m)=f(m-n,-m), \ m\neq 0.\eqno(3.46)$$ Setting $n'=m-n$,
$m'=n$ in (3.45) and using (3.46), we have
$$(-n)f(m,-m-n)+m f(m-n,-m+n)=(m-n)
f(m,-m).\eqno(3.47)$$ From (3.44) and (3.47), we have
$$2f(m,-m)-f(m-n,-m+n)=f(m+n,-m-n), \ m\neq
0.\eqno(3.48)$$ Replace $m$ by $n$ and $n$  by $m$ respectively
and use (3.43), then $$ 2f(n,-n)+f(m-n,-m+n)-2f(0,0)=f(m+n,-m-n) \
n\neq 0.\eqno(3.49)$$ From (3.47) and (3.49), we have
$$f(m+n,-m-n)-f(0,0)=f(m,-m)-f(0,0)+f(n,-n)-f(0,0).$$
We deduce that there exists $d\in\F$ such that
$$f(m,-m)=dm+f(0,0),$$
for all $m\in\z$. By (3.44), we have $$ f(m,n)=f(0,0),$$ for all
$m,n\in\z$ and $m+n\neq 0$. So

$$f(m, n)=\cases{dm+f(0,0), \quad m+n=0, \, \cr  f(0,
0),\hfill m+n\ne 0}.$$ As in the above case we can deduce that if
$f(0,0)\neq 0$, then $d=0$ and  $V\cong{\mathcal A}_{0,0,c}$,
where $c=f(0,0)\neq 0$. If $f(0,0)=0$, then
$$f(m, n)=\cases{dm, \quad m+n=0, \, \cr  0,\hfill
m+n\ne 0}.$$ Therefore $V\cong{\mathcal U}_{d}$.

\noindent{\bf Case II.}  $L(n)v_t=(t+n)v_{n+t}$ if $t\ne 0$,
$L(n)v_0=n(n+a)v_{n}$ for some $a\in\F$. We can deduce that $f(m,
t)=\cases{dm, \quad t=0 \, \cr 0,\hfill t\ne 0}$. Therefore
$V\cong{\mathcal A}(a, d)$.

\noindent{\bf Case III.}  $L(n)v_t=tv_{n+t}$ if $t\ne -n$,
$L(n)v_{-n}=-n(n+a)v_{0}$,  for some $a\in{\mathbb F}$.

We can deduce that $f(m, t)=\cases{dm, \quad t+m=0\cr\hfill 0,
\quad t+m\ne 0}$. Then $V$ is isomorphic to ${\mathcal  B}(a, d)$.

\vskip10pt

Finally suppose that $V$ is decomposable as a Vir-module. Then
$V={\mathcal A}_{0,0}'\oplus\F v_0$ by Theorem 3.3. So we can
suppose that $L(n)v_t=tv_{n+t}$ if $n+t\ne 0$, $L(n)v_{-n}=0$.
Then by calculation as in Case III we can deduce that
$I(n)v_t=\cases{dn v_0\quad  n+t=0\cr 0, \hfill n+t\ne 0}$ or
$I(n)v_t=\cases{cv_n\quad n\ne 0, t=0\cr 0,
\hfill\hbox{otherwises}}.$

In fact, in this case (3.2) becomes
$$(t+m)f(m,t)-tf(m,n+t)=m f(m+n,t),\quad m+n+t\ne 0, n+t\ne 0.\eqno(3.50)$$
$$tf(-(n+t),n+t)=(n+t)f(-t,t),\quad n+t\ne 0.\eqno(3.51)$$
$$(t+m)f(m,t)=m f(m-t,t),\quad m\ne 0.\eqno(3.52)$$

If $t\ne 0$ in (3.51), we have
$${f(-(n+t),n+t)\over n+t}={f(-t,t)\over t}.\eqno(3.53)$$
So $$f(m, -m)=dm, \quad \forall m\in\z^*$$ for some $d\in \F$.

Setting $m=-t$ in (3.52) we get
$$f(-2t, t)=0, \quad t\ne 0.$$
Setting $t=0$ in (3.51) we have $$f(0,0)=0.$$ Hence
$$f(-2t, t)=0, \quad \forall t\in\z.\eqno(3.54)$$
Setting $m=-2t, -3t, \cdots$ in (3.52) again, and by induction we
get
$$f(-kt, t)=0, \quad k\ne 1, k\in\z^+.\eqno(3.55)$$

Setting $t=0$ in (3.50) we get
$$f(m, 0)=f(n, 0), \quad m\ne 0, n\ne 0.\eqno(3.56)$$

Setting $m=-t$ in (3.50) we get
$$f(-t,n+t)=f(-t+n,t)\quad t\ne0, n+t\ne 0.\eqno(3.57)$$

Setting $n=-kt, k\ne 1, k\in\z^+$ in (3.57) we get
$$f(-t,-(k-1)t)=0,$$
it is,
$$f(t,kt)=0,\quad k\in\z^+. \eqno(3.58)$$
Using (3.55), (3.58), (3.57) and (3.52) repeatedly, we get
$$f(kt, t)=f(t, kt)=0, \quad \forall k\in \z, k\ne -1, 0.\eqno(3.59)$$

Setting $n=t$ in (3.57) and using (3.59), we get $$f(0, t)=0,
\quad t\ne 0.\eqno(3.60)$$

By $[I(m), I(n)]=\de_{m+n, 0}n C_I$, we get
$$f(n, t)f(m, n+t)=f(m, t)f(n, m+t).\eqno(3.61)$$
Setting $m=t$ in (3.61) and using (3.58), we get
$$f(n, t)f(t, n+t)=0, \quad \forall n, t\in\z. \eqno(3.62)$$

Setting $m=t\ne 0$ in (3.50)
$$f(t,n+t)=f(n+t,t)\quad  n+t\ne 0, n+2t\ne 0, t\ne 0.\eqno(3.63)$$

Setting $m=n+t$ in (3.52) we get
$$f(n+t,t)={n+t\over n+2t}f(n, t), \quad n+t\ne 0, n+2t\ne 0. \eqno(3.64)$$

Combining (3.62), (3.63) and (3.64), we get
$$f(n, t)=0, \quad  n+t\ne 0, n+2t\ne 0, t\ne 0.$$
So by (3.59) we have
$$f(n, t)=0, \quad  n+t\ne 0, t\ne 0.\eqno(3.65)$$

If $d\ne 0$, then by setting $n=-t\ne0$ in (3.62), we get
$$f(t, 0)=0.\eqno(3.66)$$
Therefore
$$f(m, n)=\cases{dm\quad  n+m=0\cr 0, \quad  n+t\ne 0}.\eqno(3.67)$$
Therefore $V\cong \tilde{\cal U}_d$.

If $d=0$, then $f(m, 0)=c$, $m\ne 0$ for some $c\in \F^*$
$$f(m, n)=\cases{c\quad  m\ne 0, n=0\cr 0, \quad  \hbox{otherwises}}.\eqno(3.68)$$
It is easy to see that $f(m,n)$ satisfies
$$
f(m,n+t)f(n,t)=f(n,m+t)f(m,t).$$  Therefore $V\cong \tilde{\cal
V}_c$ (if $c=0$, $V$ is decomposable). \hfill
$\rule[-.23ex]{1.0ex}{2.0ex}$

The above modules (if they are irreducible) or their nontrivial
simple subquotients are the Harish-Chandra modules of the
intermediate series. Moreover the above results immediately yield
the following theorem.
\begin{theo}
Any Harish-Chandra module of intermediate series over the
 Heisenberg-Virasoro algebra is isomorphic to one of
above modules (if they are irreducible) or their nontrivial simple
subquotients.
\end{theo}

\noindent{\bf Remarks.}  1. Irreducible highest weight modules of
${\cal L}$ were classified in Ref.\cite{ADKP} and Ref.\cite{B}.
Irreducible Harish-Chandra modules of ${\cal L}$ were classified
in Ref.\cite{LZ}

\begin{theo}Ref.\cite{LZ}
Any irreducible Harish-Chandra module over the twisted
Heisenberg-Virasoro algebra is a highest weight module, a lowest
weight module or a module of the intermediate series.
\end{theo}

2. The above results can be used to classify Harish-Chandra
modules of the $N = 2$ Neveu-Schwarz superalgebra, whose even part
is just the twisted Heisenberg-Virasoro algebra except for some
central elements. Especially, all Harish-Chandra modules of the
intermediate series over this superalgebra from Theorem 3.6 can be
constructed (see Ref.\cite{FJS}).

\vskip30pt \centerline{\bf ACKNOWLEDGMENTS}

\vskip15pt The part of this work was done during the visit of the
first author in the Department of Mathematics at University of
Bielefeld. He would like to expresses his special gratitude to the
`AsiaLink Project' for financial support. He is also deeply
indebted to Prof. C.M. Ringel for his kind hospitality and
continuous encouragement and instruction. The project is supported
by the NNSF (Grant 10671027, 10701019, 10571119), the ZJZSF(Grant
Y607136), and Qianjiang Excellence Project of Zhejiang Province
(No. 2007R10031).

\end{document}